\input amstex
\documentstyle{amsppt}
\output={\plainoutput}
\pageno=1
\footline={\ifnum\pageno>1 \myfl\fi}
\define\myfl{\hss\tenrm\folio\hss}
\headline={\hfil}
\NoBlackBoxes
\magnification=\magstep 1
\pagewidth{15 true cm}
\pageheight{22 true cm}

\topmatter
\title Burghelea-Friedlander-Kappeler's gluing formula for the zeta-determinant and its applications to the adiabatic decompositions of the zeta-determinant and the analytic torsion.
\endtitle
\author Yoonweon Lee
\endauthor
\affil  Department of Mathematics \\
        Inha University \\
        Inchon, 402-751, Korea
\endaffil
\endsubjclass
\thanks
The author was supported partially by Korea Research Foundation Grant KRF-2000-015-DP0045.
\endthanks
\abstract
The gluing formula of the zeta-determinant of a Laplacian given by Burghelea, Friedlander and Kappeler contains an unknown constant.
In this paper we compute this constant to complete the formula 
under the assumption of the product structure near boundary. As applications of this result,
we prove the adiabatic decomposition theorems of the zeta-determinant of a Laplacian with respect to the Dirichlet and Neumann boundary conditions
and of the analytic torsion with respect to the absolute and relative boundary conditions.
\endabstract
\endtopmatter

\document
\baselineskip 0.6 true cm plus 3 pt minus 3 pt

\S 1 {\bf Introduction}
\TagsOnRight
\vskip 0.5 true cm

In [3], Burghelea, Friedlander and Kappeler established a gluing formula for the zeta determinant
of an elliptic operator on a compact manifold. This formula contains an unknown constant which can be 
expressed in terms of the zero coefficients of some asymptotic expansions. 
In this paper we compute this constant in case that the product structure is given near boundary and 
then we apply this result to prove the adiabatic
decomposition theorems of the zeta determinant and the analytic torsion. 
Some results of this paper are known from the work of Klimek and Wojciechowski in [6] but our method is completely different from theirs.

Let $M$ be a compact oriented $m$-dimensional manifold with boundary $Z$ ($Z$ may be empty)  
and $Y$ be a hypersurface of $M$ such that
$M-Y$ has two components and $Y\cap Z=\emptyset$. We denote by $M_{1}$, $M_{2}$ the closure of each component, {\it i.e.}
$M=M_{1}\cup_{Y} M_{2}$.
Choose a collar neighborhood $N$ of $Y$, which is diffeomorphic to $[-1,1]\times Y$, $N \cap Z=\emptyset$ and choose a metric 
$g$ on $M$ which is a product metric on $N$.
Suppose that $E\rightarrow M$ is a complex vector bundle such that $E|_{N}$ has the product structure, which means that
$E|_{N}=p^{\ast}E|_{Y}$, where $p : [-1,1]\times Y\rightarrow Y$ is the canonical projection.
Let $\Delta_{M}$ be a Laplacian acting on smooth sections of $E$ and
$\Delta_{M_{1}}$, $\Delta_{M_{2}}$ be the restrictions of $\Delta_{M}$ to $M_{1}$ and $M_{2}$.
By a Laplacian we mean a positive semi-definite 2nd order differential operator whose principal symbol is 
$\sigma_{L}(\Delta_{M})(x,\xi)=||\xi||^{2}$.
We assume that  $\Delta_{M}$  is $-\partial_{u}^{2} + \Delta_{Y}$ on $N$, 
where $\partial_{u}$ is the unit normal vector
field to $Y$ on $N$, outward to $M_{1}$ and $\Delta_{Y}$ is a Laplacian on $Y$.

We denote by $D$, $B$ the Dirichlet boundary conditions on $Z$, $Y$ and 
by $C$ the Neumann boundary condition on $Y$ defined as follows.
$$D :  C^{\infty}(M_{i}) \rightarrow  C^{\infty}(Z\cap M_{i}) \qquad \text{by}\qquad
D(\phi)  =  \phi|_{Z},$$
$$B :  C^{\infty}(M_{i}) \rightarrow  C^{\infty}(Y) \qquad \text{by}\qquad
B(\phi)  =  \phi|_{Y},$$
$$\qquad C :  C^{\infty}(M_{i}) \rightarrow  C^{\infty}(Y) \qquad \text{by}\qquad
C(\phi)  =  (\partial_{u}\phi)|_{Y}.$$
Then the Laplacian $\Delta_{M,D}$ ($\Delta_{M_{i},B,D}$, $\Delta_{M_{i},C,D}$) with 
the Dirichlet condition on $Z$ (the Dirichlet condition on $Y$ and $Z$, the Neumann condition on $Y$
and the Dirichlet condition on $Z$) is 
defined by the same operator $\Delta_{M}$ ($\Delta_{M_{i}}$) with domains as follows.
$$Dom(\Delta_{M,D})= \{ \phi\in C^{\infty}(M) \mid D(\phi)=0 \},\qquad \quad $$
$$Dom(\Delta_{M_{i},B,D})= \{ \phi\in C^{\infty}(M_{i}) \mid B(\phi)=0, D(\phi)=0 \},$$
$$Dom(\Delta_{M_{i},C,D})= \{ \phi\in C^{\infty}(M_{i}) \mid C(\phi)=0, D(\phi)=0 \}.$$

For the computational reason, we consider $\Delta_{M,D_{m}}^{m}+t^{m}$, 
$\Delta_{M_{i},B_{m},D_{m}}^{m}+t^{m}$
and $\Delta_{M_{i},C_{m},D_{m}}^{m}+t^{m}$ ($t\in {\Bbb R}^{+}$) rather than 
$\Delta_{M,D}$, $\Delta_{M_{i},B,D}$
and $\Delta_{M_{i},C,D}$, where $D_{m}$, $B_{m}$ and $C_{m}$ are
the Dirichlet and the Neumann boundary conditions 
corresponding to $\Delta_{M}^{m}$, 
$\Delta_{M_{i}}^{m}$ (or $\Delta_{M}^{m}+t^{m}$, $\Delta_{M_{i}}^{m}+t^{m}$) defined as follows.
$$\align
D_{m}=(D, D\Delta_{M},\cdots, D\Delta_{M}^{m-1}), & \\
 B_{m}=(B, B\Delta_{M_{i}},\cdots, B\Delta_{M_{i}}^{m-1}), 
&\qquad C_{m}=(C, C\Delta_{M_{i}},\cdots, C\Delta_{M_{i}}^{m-1}).
\endalign
$$

Note that $$\Delta_{M,D}^{m}+t^{m}=\cases \prod_{k=-[\frac{m}{2}]}^{[\frac{m-1}{2}]}
(\Delta_{M,D}+e^{i\frac{2k\pi}{m}}t)
& \text{ if $m$ is odd} \\ 
\prod_{k=-[\frac{m}{2}]}^{[\frac{m-1}{2}]}(\Delta_{M,D}+e^{i\frac{(2k+1)\pi}{m}}t)
& \text{ if $m$ is even}.
\endcases
$$
For $-[\frac{m}{2}]\leq k \leq [\frac{m-1}{2}]$, denote $\alpha_{k}=e^{i\frac{2k\pi}{m}}$ if $m$ is odd, and
$\alpha_{k}=e^{i\frac{(2k+1)\pi}{m}}$ if $m$ is even.

Now we describe the so-called Dirichlet-to-Neumann operator $R(\alpha_{k}t):C^{\infty}(Y)\rightarrow C^{\infty}(Y)$
associated to $\Delta_{M,D}+\alpha_{k}t$ on $Y$.
Let $P_{i}(\alpha_{k}t):C^{\infty}(Y)\rightarrow C^{\infty}(M_{i})$ be the Poisson operator on $Y$ associated to 
$\Delta_{M,D}+\alpha_{k}t$, which is characterized by the following equations (for details see [3],[4],[8]).
$$ BP_{i}(\alpha_{k}t)=Id_{Y}, \qquad  DP_{i}(\alpha_{k}t)=0, \qquad
(\Delta_{M}+\alpha_{k}t) P_{i}(\alpha_{k}t)=0.$$
Then $R(\alpha_{k}t)$ is defined by the compositions of the following maps.
$$
C^{\infty}(Y) @>\delta_{ia}>> C^{\infty}(Y)\oplus C^{\infty}(Y) @>(P_{1}(\alpha_{k}t), P_{2}(\alpha_{k}t))>>
C^{\infty}(M_{1})\oplus C^{\infty}(M_{2}) $$
$$@>(C_{1},C_{2})>> C^{\infty}(Y)\oplus C^{\infty}(Y)
@>\delta_{if}>> C^{\infty}(Y),$$
where $\delta_{ia}(g)=(g,g)$, $C_{1}(\phi_{1})=(\partial_{u}\phi_{1})|_{Y}$, $C_{2}(\phi_{2})=(\partial_{u}\phi_{2})|_{Y}$
and $\delta_{if}(g,h)=g-h$.
It is known that $R(\alpha_{k}t)$ is a $\Psi$DO of order $1$ ({\it c.f.} Theorem 2.1)
and by choosing $\pi$ as an Agmon angle, 
$\log DetR(\alpha_{k}t)$ is well defined.
The following theorem is due to Burghelea, Friedlander and Kappeler ([8], see also [3] and [4]).

\proclaim{Theorem 1.1}
$$\log Det(\Delta_{M,{D_{m}}}^{m}+t^{m}) - \log Det(\Delta_{M_{1},B_{m},D_{m}}^{m} +t^{m}) -  
\log Det(\Delta_{M_{2},B_{m},D_{m}}^{m}+t^{m})$$
$$=-\sum_{k=-[\frac{m}{2}]}^{[\frac{m-1}{2}]}c_{k}+
\sum_{k=-[\frac{m}{2}]}^{[\frac{m-1}{2}]}\log DetR(\alpha_{k}t),$$
where $c_{k}$ is the zero coefficient in the asymptotic expansion of $\log DetR(\alpha_{k}t)$ as $t\rightarrow\infty$.
\endproclaim

{\it Remark} : \hskip 0.2 true cm
In [3] and [8], Theorem 1.1 was proved only in case that $Z=\emptyset$. However, the proof can be extended 
without any modification to the case that $Z$ is non-empty. 
\vskip 0.2 true cm

The purpose of this paper is to compute the zero coefficients in Theorem 1.1 under the assumption of the product
structures on $N$, $E|_{N}$ and then
we apply this result to prove the adiabatic decomposition theorems 
of the zeta-determinant of a Laplacian and the analytic torsion. 
We first have the following theorem.
\proclaim{Theorem 1.2}
We assume the product structures of $M$ and $E$ on $N$ and $\Delta_{M}=-\partial_{u}^{2}+\Delta_{Y}$ on $N$. Then :

\vskip 0.2 true cm

$\sum_{k}c_{k} =  m \log2 \cdot (\zeta_{\Delta_{Y}}(0)+dimker\Delta_{Y}).$

\endproclaim
\vskip 0.2 true cm
Setting $t=0$, we get the following corollary.
\proclaim{Corollary 1.3}
We further assume that $\Delta_{M,D}$ is invertible. 
Then :
\vskip 0.1 true cm
 $\log Det\Delta_{M,D} - \log Det\Delta_{M_{1},B,D} - \log Det\Delta_{M_{2},B,D} =$

\vskip 0.2 true cm
 
$\qquad \qquad \qquad \qquad \qquad \qquad 
- \log 2 \cdot (\zeta_{\Delta_{Y}}(0) +dimker\Delta_{Y})+\log Det R$.
\endproclaim

\vskip 0.2true cm
{\it Remark} : \hskip 0.2 true cm
(1) \hskip 0.1 true cm 
If $dimY$ is odd, it is well-known that $\zeta_{\Delta_{Y}}(0)+dimker\Delta_{Y}=0$. In this case,
the assertion in Corollary 1.3 can be written as follows.
 $$\log Det\Delta_{M,D} - \log Det\Delta_{M_{1},B,D} - \log Det\Delta_{M_{2},B,D} = \log Det R,$$
which was observed in [7].  \newline
(2) \hskip 0.1 true cm
Theorem 1.1, Theorem 1.2 and Corollary 1.3 also hold when we impose the
absolute (or the relative) boundary condition on $Z$ (see Theorem 5.2).

\vskip 0.2 true cm

The main idea of proving Theorem 1.2 is to show that under the assumption of the product structure,
$R(\alpha_{k}t)$ can be expressed as $2 \sqrt{\Delta_{Y}+\alpha_{k}t}$ $ + $ a smoothing operator
(Theorem 2.1). 
We are going to show this fact in the next section
by using an observation due to I.M. Gelfand, (probably unpublished), that the Dirichlet-to-Neumann operator
satisfies a Ricatti type equation ({\it c.f.} (2.2)). 

\vskip 0.2 true cm

Now we apply Corollary 1.3 to discuss the adiabatic decomposition of the zeta-determinant of a Laplacian 
into the zeta-determinants of Laplacians with the Dirichlet
and Neumann boundary conditions. 
Recall that $N$ is a collar neighborhood of $Y$, which is diffeomorphic to $[-1,1]\times Y$. We denote by $M_{r}$ the compact manifold
with boundary obtained by attaching $N_{r+1}=[-r-1,r+1]\times Y$ on 
$M-(-\frac{1}{2},\frac{1}{2})\times Y$ by identifying $[-1,-\frac{1}{2}]\times Y$ with 
$[-r-1,-r-\frac{1}{2}]\times Y$
and $[\frac{1}{2},1]\times Y$ with $[r+\frac{1}{2},r+1]\times Y$. 
We also denote by $M_{1,r}$, $M_{2,r}$ the manifolds with boundary which
are obtained by attaching $[-r,0]\times Y$, $[0,r]\times Y$ on $M_{1}$, $M_{2}$ by identifying $Y$ with $\{ -r \}\times Y$
and $Y$ with $\{ r \}\times Y$, respectively. 
Then the bundle $E\rightarrow M$ and the Laplacian $\Delta_{M}$ on $M$ can be extended naturally to 
the bundle $E_{r}\rightarrow M_{r}$ and the Laplacian $\Delta_{M_{r}}$ on $M_{r}$.  

To describe the next result, we need to define the operators $Q_{i}:C^{\infty}(Y) \rightarrow C^{\infty}(Y)$ ($i=1$, $2$) by 
slightly modifying the Dirichlet-to-Neumann operator.
For $f\in C^{\infty}(Y)$, choose $\phi_{i}\in C^{\infty}(M_{i})$ satisfying $\Delta_{M_{i}}\phi_{i}=0$, $\phi_{i}|_{Z}=0$ and
$\phi_{i}|_{Y}=f$. We define 
$$Q_{1}(f)=(\partial_{u}\phi_{1})|_{Y} , \qquad \qquad Q_{2}(f)=(-\partial_{u}\phi_{2})|_{Y}.$$
Then each $Q_{i}$ is an elliptic $\Psi DO$ of order 1 ({\it c.f.} Theorem 2.1) and
the Dirichlet-to-Neumann operator $R$ is $R=Q_{1}+Q_{2}$.  
The following is the second result of this paper.

\proclaim{Theorem 1.4}
We assume that both $Q_{1}+\sqrt{\Delta_{Y}}$ and $Q_{2}+\sqrt{\Delta_{Y}}$ are invertible operators and
$k=dimKer\Delta_{Y}$. We further assume that $\Delta_{M_{r},D}$ is invertible for $r$ large enough.
Then :
$$\multline
\lim_{r\to\infty}\left\{ \log Det(\Delta_{M_{r},D})-\log Det(\Delta_{M_{1,r},B,D})-\log Det(\Delta_{M_{2,r},B,D}) 
+k \log r\right\} \\
= \frac{1}{2}\log Det\Delta_{Y}.\endmultline 
$$
\endproclaim

{\it Remark} : \hskip 0.2 true cm
(1) If $\Delta_{Y}$ has non-trivial kernel, we define $Det\Delta_{Y}$ 
from the zeta function $\zeta_{\Delta_{Y}}(s)$ consisting of only non-zero eigenvalues. \newline
(2) If $\Delta_{M}$ is a connection Laplacian for a connection compatible with the inner product, each $Q_{i}$ is a
non-negative operator (Lemma 4.3). \newline
(3) Suppose that $\Delta_{M}=A^{2}$ for a Dirac operator $A$ which has the form $G(\partial_{u}+B)$ near $Y$ with $G$ 
a bundle automorphism satisfying
$$G^{\ast}=-G, \qquad G^{2}=-Id, \qquad B^{\ast}=B,  \qquad GB=-BG.   \tag1.1$$
Here $G$ and $B$ do not depend on the normal coordinate $u$.
Then the invertiblity of both $Q_{1}+\sqrt{B^{2}}$ and $Q_{2}+\sqrt{B^{2}}$ is equivalent to the non-existence of the extended 
$L^{2}$-solutions  of $A_{M_{1,\infty}}$, $A_{M_{2,\infty}}$
on ${M_{1,\infty}}$ and ${M_{2,\infty}}$ (Corollary 4.5). \newline 
(4) Suppose that $\Delta_{M}$ is a connection Laplacian or a Dirac Laplacian for a connection compatible with the inner product
and $\Delta_{M,D}$ is invertible. Then the invertiblity of both $Q_{1}+\sqrt{\Delta_{Y}}$ and $Q_{2}+\sqrt{\Delta_{Y}}$
implies the invertiblity of $\Delta_{M_{r},D}$ for $r$ large enough (Lemma 4.6).

\vskip 0.3 true cm

Let ${\tilde M_{1,r}}$ be the double of $M_{1,r}$. Then it is a well-known fact that
$$\log Det\Delta_{{\tilde M_{1,r}},D,D}=\log Det\Delta_{M_{1,r},C,D}+\log Det\Delta_{M_{1,r},B,D}.$$
Combining this fact with Corollary 1.3 and Theorem 1.4, we have the following result.

\proclaim{Corollary 1.5} 
\hskip 0.2 true cm
We assume the hypotheses in Theorem 1.4. Then :
$$
(1) \lim_{r\to\infty}\left\{ \log Det(\Delta_{M_{1,r},C,D})-\log Det(\Delta_{M_{1,r},B,D}) + k \log r \right\} 
 = \frac{1}{2}\log Det(\Delta_{Y}).
$$
$$
(2) \lim_{r\to\infty}\left\{ \log Det(\Delta_{M_{r},D})-\log Det(\Delta_{M_{1,r},C,D})-\log Det(\Delta_{M_{2,r},B,D}) \right\}=0.
\qquad $$
\endproclaim

\vskip 0.3 true cm

Finally we discuss the adiabatic decomposition of the analytic torsion
into the analytic torsions with the absolute and relative boundary
conditions.

Here we assume that $M$ is a closed manifold with a hypersurface $Y$ and $M$ has a product structure
near $Y$. We define $M_{r}$, $M_{1,r}$ and $M_{2,r}$ as above so that 
$M_{r}=M_{1,r} \cup_{\{ 0 \}\times Y} M_{2,r}$. 
Suppose that $\rho_{M_{r}}$ ($\rho_{M_{1,r}}$, $\rho_{M_{2,r}}$, $\rho_{Y}$) is an orthogonal representation
of $\pi_{1}(M_{r})$ ($\pi_{1}(M_{1,r})$, $\pi_{1}(M_{2,r})$, $\pi_{1}(Y)$) to $SO(n)$, respectively.
Then we can define the analytic torsions $\tau(M_{r},\rho_{M_{r}})$, $\tau_{abs}(M_{i,r},\rho_{M_{i,r}})$, 
$\tau_{rel}(M_{i,r},\rho_{M_{i,r}})$ 
$(i=1,2)$, $\tau(Y,\rho_{Y})$ in the standard way (for the definitions, see Section 5). 
Our goal is to recover Klimek-Wojciechowski's result about the analytic torsion in [6] as follows.

First, let us consider $M_{1,r}$ (a manifold with boundary $Y$) only. For a given representation
$\rho_{M_{1,r}} : \pi_{1}(M_{1,r}) \rightarrow SO(n)$ and the natural homomorphism 
$\iota_{Y} : \pi_{1}(Y) \rightarrow \pi_{1}(M_{1,r})$, define $\rho_{Y} : \pi_{1}(Y) \rightarrow SO(n)$
by $\rho_{Y}=\rho_{M_{1,r}} \circ \iota_{Y}$. We denote by $\Delta_{Y}^{q}$ ($\Delta_{M_{1,r}}^{q}$)
the Hodge Laplacian acting on $q$-forms on $Y$ (on $M_{1,r}$) valued in $E_{\rho_{Y}}$ ($E_{\rho_{M_{1,r}}}$),
where $E_{\rho_{Y}}= {\tilde Y} \times_{\rho_{Y}} {\Bbb R}^{n}$ with ${\tilde Y}$ the universal covering space of $Y$
($E_{\rho_{M_{1,r}}}$ is defined in the same way).
We define $Q_{1}^{q}$ by the same way as in Theorem 1.4 with the bundle 
$E=\wedge^{q}T^{\ast}M_{1,r}\otimes E_{\rho_{M_{1,r}}}$. 
If necessary, by tensoring ${\Bbb C}$ on $E$, we regard $E$ as a complex vector bundle.
Then we have the following theorem.

\proclaim{Theorem 1.6}
Suppose that for each $q$, 
$Q_{1}^{q} +\left( \smallmatrix \sqrt{\Delta_{Y}^{q}} & 0 \\  0 & \sqrt{\Delta_{Y}^{q-1}} \endsmallmatrix \right)$ is 
an invertible operator on $\{ -r \}\times Y$
and $H^{q}(M_{1,r};\rho_{M_{1,r}})$, $H^{q}(M_{1,r}, Y ; \rho_{M_{1,r}})$ are trivial groups.
Then,
$$ 
\lim_{r \to \infty} \left\{ \log\tau_{abs}(M_{1,r},\rho_{M_{1,r}}) - \log\tau_{rel}(M_{1,r},\rho_{M_{1,r}})  \right\}
= \log\tau(Y;\rho_{Y}). $$
\endproclaim
\vskip 0.2 true cm
{\it Remark} : \hskip 0.3 true cm
If $Q_{1}^{q} +\left( \smallmatrix \sqrt{\Delta_{Y}^{q}} & 0 \\  0 & \sqrt{\Delta_{Y}^{q-1}} \endsmallmatrix \right)$ is 
invertible, by Corollary 4.5 there is no extended $L^{2}$-solutions of $d_{q}+d_{q}^{\ast}$ on $M_{1,\infty}$, which implies that
$Ker\Delta^{q-1}_{Y}=Ker\Delta^{q}_{Y}=0$ ({\it c.f.} [1], [2], [5]).

\vskip 0.3 true cm
Next, we consider the closed manifold $M_{r}$ and manifolds with boundary $M_{i,r}$, $(i=1,2)$. For a given representation
$\rho_{M_{r}} : \pi_{1}(M_{r}) \rightarrow SO(n)$ and the natural homomorphisms
$\iota_{M_{i,r}} : \pi_{1}(M_{i,r}) \rightarrow \pi_{1}(M_{r})$, 
$\iota_{Y} : \pi_{1}(Y) \rightarrow \pi_{1}(M_{i,r})$,
define 
$\rho_{M_{i,r}} : \pi_{1}(M_{i,r}) \rightarrow SO(n)$,
$\rho_{Y} : \pi_{1}(Y) \rightarrow SO(n)$,
by $\rho_{M_{i,r}}=\rho_{M_{r}}\circ \iota_{M_{i,r}}$, $\rho_{Y}=\rho_{M_{i,r}}\circ \iota_{Y}$.
We also define $\Delta_{Y}^{q}$, $Q_{1}^{q}$ and $Q_{2}^{q}$ as in Theorem 1.6.

\proclaim{Theorem 1.7}
Suppose that for each $q$, $Q_{1}^{q} +\left( \smallmatrix \sqrt{\Delta_{Y}^{q}} & 0 \\  0 & \sqrt{\Delta_{Y}^{q-1}} \endsmallmatrix \right)$ 
and  \newline
$Q_{2}^{q} +\left( \smallmatrix \sqrt{\Delta_{Y}^{q}} & 0 \\  0 & \sqrt{\Delta_{Y}^{q-1}} \endsmallmatrix \right)$ are invertible operators 
on $\{ -r \}\times Y$, $\{ r \}\times Y$, and  \newline 
$H^{q}(M_{r};\rho_{M_{r}})$, $H^{q}(M_{1,r};\rho_{M_{1,r}})$, $H^{q}(M_{2,r},Y;\rho_{M_{2,r}})$ are trivial groups.
Then :
$$ (1) \lim_{r\to\infty}\left( \log Det\Delta_{M_{r}}^{q}-\log Det\Delta^{q}_{M_{1,r},abs}-\log Det\Delta^{q}_{M_{2,r},rel} \right) = 0. \qquad 
\qquad  $$

$$ (2) \lim_{r\to\infty}\left(\log \tau(M_{r};\rho_{M_{r}}) -\log \tau_{abs}(M_{1,r};\rho_{M_{1,r}})
- \log\tau_{rel}(M_{2,r};\rho_{M_{2,r}}) \right) = 0.$$
\endproclaim

\vskip 0.3 true cm

{\it Remark} : \hskip 0.2 true cm
Recently J. Park and K. Wojciechowski showed the following result in [10]. Suppose that $M$ is an odd-dimensional compact manifold
with $M=M_{1}\cup_{Y}M_{2}$ and $D$ is a Dirac operator acting on smooth sections of a Clifford module bundle $E$
with $D=G(\partial_{u}+B)$ near $Y$. 
Denote by $P_{>}$, $P_{<}$
the Atiyah-Patodi-Singer boundary conditions projecting the positive and negative eigenspaces of $B$, respectively.
Assume that 
$$KerB =\{ 0 \}, \qquad \qquad Ker_{L^{2}}D_{1,\infty}=Ker_{L^{2}}D_{2,\infty}= \{ 0 \},$$
where $Ker_{L^{2}}D_{i,\infty}$ is the set of all extended $L^{2}$-solutions of $D_{i,\infty}$ on $M_{i,\infty}$.
Then
$$\lim_{r\to\infty}\left\{ \log DetD_{r}^{2}-\log DetD^{2}_{M_{1,r},P_{>}}-\log DetD^{2}_{M_{2,r},P_{<}} \right\}
= - \log 2  \cdot \zeta_{B^{2}}(0).$$
This result is the main motivation of this paper. In the next work ([9]) we are going to recover this result by using the techniques in this paper.

\vskip 1 true cm

\S 2 {\bf Asymptotic symbol of  $R(\alpha_{k}t)$ }

\vskip 0.5 true cm

In this section, we are going to describe the asymptotic symbol of $R(\alpha_{k}t)$.
The following method is observed by I.M. Gelfand.

We start from defining $Q_{i}(\alpha_{k}t) : C^{\infty}(Y)\rightarrow C^{\infty}(Y)$ ($i=1, 2$) as follows.
For $f\in C^{\infty}(Y)$, choose $\phi_{i}\in C^{\infty}(M_{i})$ such that
$$(\Delta_{M_{i}}+\alpha_{k}t)\phi_{i}=0, \quad \phi_{i}|_{Y}=f, \quad \phi_{i}|_{Z}=0.$$ 
Then we define
$$Q_{1}(\alpha_{k}t)(f)=(\partial_{u}\phi_{1})|_{Y} \text{ and } 
Q_{2}(\alpha_{k}t)(f)=(-\partial_{u}\phi_{1})|_{Y}.$$
From this definition, we get
$$R(\alpha_{k}t)=Q_{1}(\alpha_{k}t) + Q_{2}(\alpha_{k}t)$$
and it's enough to consider $Q_{1}(\alpha_{k}t)$ only. From now on we denote $Q_{1}(\alpha_{k}t)$ simply by
$Q(\alpha_{k}t)$.

For $f\in C^{\infty}(Y)$, let $\varphi$ be a solution of $\Delta_{M_{1}}+\alpha_{k}t$ with  
$\varphi|_{Y}=f$ and $\phi|_{Z}=0$. Then,
$$\frac{d}{du}\varphi(u,y) = Q_{u}(\alpha_{k}t) \varphi(u,y),$$
where $Q_{u}(\alpha_{k}t)$ is defined similarly as $Q(\alpha_{k}t) = Q_{0}(\alpha_{k}t)$ 
at the level $\{ u \}\times Y$.
$$
\frac{d^{2}}{du^{2}}\varphi(u,y)  = \left(\frac{d}{du} Q_{u}(\alpha_{k}t)\right) \varphi(u,y) + 
Q_{u}(\alpha_{k}t)^{2} \varphi(u,y).$$
For $0\leq u < 1$,
$$
(\Delta_{Y}+ \alpha_{k}t) \varphi(u,y) = \left(\frac{d}{du} Q_{u}(\alpha_{k}t)\right) \varphi(u,y) + 
Q_{u}(\alpha_{k}t)^{2} \varphi(u,y).$$
Consequently for $0\leq u < 1$,
$$
\frac{d}{du} Q_{u}(\alpha_{k}t) = - Q_{u}(\alpha_{k}t)^{2} + (\Delta_{Y}+ \alpha_{k}t).  \tag2.2 
$$

Now let us consider the asymptotic symbol of $Q_{u}(\alpha_{k}t)$ as follows.
$$\sigma(Q_{u}(\alpha_{k}t)) \sim q_{1}(u,y,\xi) + q_{0}(u,y,\xi)
+ \cdots + q_{1-j}(u,y,\xi) + \cdots, $$
where $q_{1-j}(u,y,\xi)$ is the homogeneous part of $\sigma(Q_{u}(\alpha_{k}t))$ of order $1-j$ 
with respect to $\xi$. Then
$$
\sigma\left(\frac{d}{du}Q_{u}(\alpha_{k}t)\right) \sim \frac{d}{du} q_{1}(u,y,\xi) + 
\frac{d}{du} q_{0}(u,y,\xi)
+ \cdots + \frac{d}{du}q_{1-j}(u,y,\xi) + \cdots.  \tag2.3 
$$
Note that
$$
\split
\sigma\left(Q_{u}(\alpha_{k}t)^{2}\right) &  \sim \sum_{k=0}^{\infty}\sum \Sb |\omega|+i+j=k \\i, j\geq 0
\endSb \frac{1}{\omega !}d^{\omega}_{\xi}q_{1-i}(u,y,\xi) \cdot
D_{y}^{\omega}q_{1-j}(u,y,\xi) \\
& = q_{1}^{2}(u,y,\xi) + (d_{\xi}q_{1} D_{y}q_{1} 
+q_{0} q_{1} + q_{1} q_{0}) + \cdots .  
\endsplit   \tag2.4
$$
Suppose that
$$
\sigma(\Delta_{Y}+\alpha_{k}t) = (p_{2}(y,\xi)+\alpha_{k}t Id) + 
p_{1}(y,\xi) + p_{0}(y,\xi). 
$$
Since $\frac{d}{du}Q_{u}(\alpha_{k}t)$ is a $\Psi$DO of order $1$, 
$q_{1}^{2}(u,y,\xi)=p_{2}(y,\xi)+\alpha_{k}t Id$.
Applying the argument of Lemma 3.3 in [8] to the double of a manifold with boundary, one can show that 
$$q_{1}(u,y,\xi)= \sqrt{p_{2}(y,\xi)+\alpha_{k}t Id}.  \tag2.5 $$
Hence $q_{1}$ does not depend on $u$ and $\frac{d}{du}q_{1}(u,y,\xi)=0$.
Again, from (2.2), (2.3) and (2.4), since $q_{1}$ is a scalar matrix,
$(d_{\xi}q_{1} D_{y}q_{1} + 2 q_{1} q_{0}) = 
p_{1}(y, \xi)$ and
$$q_{0}(u,y,\xi) = (2q_{1}(y,\xi))^{-1} \left( p_{1}(y, \xi) - 
d_{\xi}q_{1}(y,\xi)\cdot  D_{y}q_{1}(y,\xi) \right) .$$
Hence $q_{0}(u,y,\xi)$ does not depend on $u$
and $\frac{d}{du}q_{0}(u,y,\xi)=0$. In general,
$$
q_{-1}= (2q_{1})^{-1}\left\{ -\sum\Sb |\omega|+i+j=2 \\ 0\leq i, j\leq 1 \endSb
\frac{1}{\omega !}d^{\omega}_{\xi}q_{1-i}(y,\xi) \cdot
D_{y}^{\omega}q_{1-j}(y,\xi) + p_{0}(y,\xi) \right\}$$
and for $k\geq 3$,
$$
q_{1-k}= (2q_{1})^{-1}\left\{ -\sum\Sb |\omega|+i+j=k \\ 0\leq i, j\leq k-1 \endSb
\frac{1}{\omega !}d^{\omega}_{\xi}q_{1-i}(y,\xi) \cdot
D_{y}^{\omega}q_{1-j}(y,\xi) \right\}.$$
Hence, each $q_{1-k}$ does not depend on $u$ and this implies that
$\frac{d}{du}Q_{u}(\alpha_{k}t)$ is a smoothing
operator. Setting $u=0$ in (2.2),
$$Q(\alpha_{k}t)^{2} = (\Delta_{Y}+\alpha_{k}t) + \text{ a smoothing operator } \tag2.6 $$
and we get the following theorem.
\proclaim{Theorem 2.1}
Under the assumption of the product structure near $N$, we have the followings.
\vskip 0.1 true cm
(1) $ Q(\alpha_{k}t) = \sqrt{\Delta_{Y}+\alpha_{k}t} + \text{ a smoothing operator }. $ \newline
\vskip 0.1 true cm
(2) $ R(\alpha_{k}t) = 2 \sqrt{\Delta_{Y}+\alpha_{k}t} + \text{ a smoothing operator }. $
\endproclaim
{\it Proof} 
\hskip 0.3 true cm
It's enough to show the first statement.
From (2.5) we have
$$Q(\alpha_{k}t) = \sqrt{\Delta_{Y}+\alpha_{k}t} + A,$$
where $A$ is an operator of order 0. Squaring both sides and using (2.6), we have
$$
\split
Q(\alpha_{k}t)^{2} & = (\Delta_{Y}+\alpha_{k}t) + 
\sqrt{\Delta_{Y}+\alpha_{k}t}A + A\sqrt{\Delta_{Y}+\alpha_{k}t} + A^{2}  \\ 
 & = (\Delta_{Y}+\alpha_{k}t ) + \text{ a smoothing operator }.
\endsplit
$$
Hence $\sqrt{\Delta_{Y}+\alpha_{k}t}A + A\sqrt{\Delta_{Y}+\alpha_{k}t} + A^{2}$ is a smoothing operator,
which implies that $A$ is a smoothing operator.
\qed

\vskip 1 true cm
\S 3 {\bf The computation of the zero coefficient of $\log DetR(\alpha_{k}t)$ as $t\rightarrow\infty$ }
\vskip 0.5 true cm
\rm

It is shown in [3] that $\log DetR(\alpha_{k}t)$ has an asymptotic expansion as $t\rightarrow\infty$
and each coefficient can be computed by the asymptotic symbol of $R(\alpha_{k}t)$. Hence, from Theorem 2.1,
$\log DetR(\alpha_{k}t)$ and $\log Det(2\sqrt{\Delta_{Y}+\alpha_{k}t})$ have the same asymptotic expansions as
$t\rightarrow\infty$. In this section, we are going to compute the asymptotic expansion of 
$\log Det(2\sqrt{\Delta_{Y}+\alpha_{k}t})$ by using the method in [12].

Note that 
$$\log Det(2\sqrt{\Delta_{Y}+\alpha_{k}t}) = \log2\cdot\zeta_{(\Delta_{Y}+\alpha_{k}t)}(0) +
 \frac{1}{2} \log Det(\Delta_{Y}+\alpha_{k}t)   \tag3.1$$
and we are going to consider $\log Det(\Delta_{Y}+\alpha_{k}t)$.
Since $Re(\alpha_{k})$ is possibly negative, we avoid this difficulty as follows.
Put  $\alpha_{k}=e^{i\theta_{k}}$ with $\theta_{k}=\frac{2k\pi}{m}$ for $m$ odd and
$\frac{(2k+1)\pi}{m}$ for $m$ even. 
Choose an angle $\phi_{k}$ with $0\leq |\phi_{k}| < \frac{\pi}{2}$ so that 
$Re (e^{i(\theta_{k}-\phi_{k})}) > 0$. 
(In fact, if $0 \leq |\theta_{k}| <\frac{\pi}{2}$, we choose $\phi_{k}=0$.) 
Then
$$
\aligned
\log Det(\Delta_{Y}+\alpha_{k}t) = \log Det \{ e^{i\phi_{k}}(e^{-i\phi_{k}}\Delta_{Y}+
e^{i(\theta_{k}-\phi_{k})}t)\}  \\
=-\frac{d}{ds}|_{s=0}\left\{e^{-i\phi_{k}s}\zeta_{(e^{-i\phi_{k}}\Delta_{Y}+
e^{i(\theta_{k}-\phi_{k})}t)}(s) \right\}  \\
=i\phi_{k}\zeta_{(e^{-i\phi_{k}}\Delta_{Y}+e^{i(\theta_{k}-\phi_{k})}t)}(0) + 
\log Det(e^{-i\phi_{k}}\Delta_{Y}+e^{i(\theta_{k}-\phi_{k})}t).   
\endaligned   \tag3.2
$$
Put ${\tilde \theta_{k}}=\theta_{k}-\phi_{k}$. Then
$$
\aligned
\zeta_{(e^{-i\phi_{k}}\Delta_{Y}+e^{i(\theta_{k}-\phi_{k})}t)}(s)
=\frac{1}{\Gamma (s)}\int_{0}^{\infty}r^{s-1}
Tre^{-r(e^{-i\phi_{k}}\Delta_{Y}+e^{i{\tilde \theta_{k}}}t)}dr  \\
=\frac{1}{\Gamma (s)}\int_{0}^{\infty}r^{s-1}
e^{-rte^{i{\tilde \theta_{k}}}} Tre^{-re^{-i\phi_{k}}\Delta_{Y}}dr.
\endaligned
$$
The following lemma is a well-known fact.

\proclaim{Lemma 3.1}
As $r\rightarrow 0$, we have the following asymptotic expansion
$$
Tre^{-re^{-i\phi_{k}}\Delta_{Y}}\sim b_{1}r^{-\frac{m-1}{2}} + b_{2}r^{-\frac{m-2}{2}} + \cdots +
b_{m} + b_{m-1}r^{\frac{1}{2}} + \cdots  $$
with $b_{m}=\zeta_{\Delta_{Y}}(0) + dimKer\Delta_{Y}$.
\endproclaim

Now we are going to compute the asymptotic expansion of $\zeta_{(e^{-i\phi_{k}}\Delta_{Y}+e^{i{\tilde \theta_{k}}}t)}(s)$
as $t\rightarrow\infty$.
$$
\aligned
\zeta_{(e^{-i\phi_{k}}\Delta_{Y}+e^{i{\tilde \theta_{k}}}t)}(s)
= \frac{1}{\Gamma (s)}\int_{0}^{\infty}r^{s-1}
e^{-rte^{i{\tilde \theta_{k}}}} Tre^{-re^{-i\phi_{k}}\Delta_{Y}}dr \\
= \frac{1}{\Gamma (s)}\int_{0}^{\infty}(\frac{u}{t})^{s-1}
e^{-ue^{i{\tilde \theta_{k}}}} Tre^{-\frac{u}{t}e^{-i\phi_{k}}\Delta_{Y}}\frac{1}{t}du \\
= t^{-s} \frac{1}{\Gamma (s)}\int_{0}^{\infty}u^{s-1}
e^{-ue^{i{\tilde \theta_{k}}}} Tre^{-\frac{u}{t}e^{-i\phi_{k}}\Delta_{Y}}du.
\endaligned
$$
As $t\rightarrow\infty$,
$$
\aligned
\zeta_{(e^{-i\phi_{k}}\Delta_{Y}+e^{i{\tilde \theta_{k}}}t)}(s)
 \sim t^{-s} \sum_{j=1}^{\infty}\frac{1}{\Gamma (s)}b_{j}\int_{0}^{\infty}u^{s-1}
(\frac{u}{t})^{\frac{j-m}{2}}e^{-ue^{i{\tilde \theta_{k}}}} du \\
 = \sum_{j=1}^{\infty}b_{j}t^{-s+\frac{m-j}{2}}\frac{1}{\Gamma (s)}\int_{0}^{\infty}u^{s+\frac{j-m}{2}-1}
e^{-ue^{i{\tilde \theta_{k}}}} du \\
 =  \sum_{j=1}^{\infty}b_{j}t^{-s+\frac{m-j}{2}}\frac{1}{\Gamma (s)}(e^{-i{\tilde \theta_{k}}})^{s+\frac{j-m}{2}}
\int_{0}^{\infty} (ue^{i{\tilde \theta_{k}}})^{s+\frac{j-m}{2}-1}
e^{-ue^{i{\tilde \theta_{k}}}}(e^{i{\tilde \theta_{k}}}) du.
\endaligned
$$  

Consider the contour integral $\int_{C}z^{s+\frac{j-m}{2}-1}e^{-z}dz$ for $Re s > \frac{m-j}{2}$, where
$$ \multline
C= \{ re^{i{\tilde \theta_{k}}} \mid \epsilon \leq r\leq R \}  \cup  \{ \epsilon e^{i\theta} \mid 0\leq \theta\leq{\tilde \theta_{k}} \} \\
\cup \{ r \mid \epsilon \leq r\leq R \} \cup  \{ Re^{i\theta} \mid 0\leq \theta\leq{\tilde \theta_{k}} \}
\endmultline $$
and oriented counterclockwise. Then one can check that
$$
\int_{0}^{\infty} (ue^{i{\tilde \theta_{k}}})^{s+\frac{j-m}{2}-1}
e^{-ue^{i{\tilde \theta_{k}}}}(e^{i{\tilde \theta_{k}}}) du = \int_{0}^{\infty}r^{s+\frac{j-m}{2}-1}e^{-r}dr
=\Gamma (s+\frac{j-m}{2}).$$
We, therefore, obtain the following asymptotic expansion for $t\rightarrow\infty$.
$$
\aligned
\zeta_{(e^{-i\phi_{k}}\Delta_{Y}+e^{i{\tilde \theta_{k}}}t)}(s)
 \sim \sum_{j=1}^{\infty}b_{j}(e^{-i{\tilde \theta_{k}}})^{s+\frac{j-m}{2}}
\frac{\Gamma(s+\frac{j-m}{2})}{\Gamma(s)}t^{-s+\frac{m-j}{2}} \\
= s \sum^{\infty}\Sb j=1 \\ j\neq m\endSb b_{j}(e^{-i{\tilde \theta_{k}}})^{s+\frac{j-m}{2}}
\frac{\Gamma(s+\frac{j-m}{2})}{\Gamma(s+1)}t^{-s+\frac{m-j}{2}} + 
b_{m}e^{-i{\tilde \theta_{k}}s}t^{-s}.
\endaligned
$$
This gives the asymptotic expansion of $\zeta_{(e^{-i\phi_{k}}\Delta_{Y}+e^{i{\tilde \theta_{k}}}t)}(s)$ as $t\rightarrow\infty$.
In view of Theorem 1.1 we are mainly interested in the zero coefficients in the asymptotic expansions of 
$\zeta_{(e^{-i\phi_{k}}\Delta_{Y}+e^{i{\tilde \theta_{k}}}t)}(0)$ and 
$\zeta_{(e^{-i\phi_{k}}\Delta_{Y}+e^{i{\tilde \theta_{k}}}t)}^{\prime}(0)$ as $t\rightarrow\infty$.

First, setting $s=0$, the zero coefficient $\pi_{0}(\zeta_{(e^{-i\phi_{k}}\Delta_{Y}+e^{i{\tilde \theta_{k}}}t)}(0))$ in the asymptotic
expansion of $\zeta_{(e^{-i\phi_{k}}\Delta_{Y}+e^{i{\tilde \theta_{k}}}t)}(0)$ is the following.
$$\pi_{0}(\zeta_{(e^{-i\phi_{k}}\Delta_{Y}+e^{i{\tilde \theta_{k}}}t)}(0))= b_{m}=
\zeta_{\Delta_{Y}}(0)+dimker\Delta_{Y}.  \tag3.3$$
Taking derivative at $s=0$, the zero coefficient of $\zeta_{(e^{-i\phi_{k}}\Delta_{Y}+e^{i{\tilde \theta_{k}}}t)}^{\prime}(0)$
can be obtained only in the term $b_{m}e^{-i{\tilde \theta_{k}}s}t^{-s}$. 
Hence, by (3.2) and (3.3), the zero coefficient $\pi_{0}(\Delta_{Y}+\alpha_{k}t)$ 
in the asymptotic expansion of $\log Det(\Delta_{Y}+\alpha_{k}t)$
as $t\rightarrow\infty$ is the following.
$$
\split
\pi_{0}(\Delta_{Y}+\alpha_{k}t)
& = i\phi_{k}(\zeta_{\Delta_{Y}}(0)+dimker\Delta_{Y})+i(\theta_{k}-\phi_{k})
(\zeta_{\Delta_{Y}}(0)+dimker\Delta_{Y}) \\
& = i\theta_{k}(\zeta_{\Delta_{Y}}(0)+dimker\Delta_{Y}).  
\endsplit   \tag3.4
$$
We summarize the above computations as follows.
\proclaim{Proposition 3.2}
The zero coefficients in the asymptotic expansions of      \newline
$\zeta_{(e^{-i\phi_{k}}\Delta_{Y}+e^{i{\tilde \theta_{k}}}t)}(0)$
and of $\log Det(\Delta_{Y}+\alpha_{k}t)$ as $t\rightarrow\infty$ are the followings.

\vskip 0.2 true cm

(1) $\pi_{0}(\zeta_{(e^{-i\phi_{k}}\Delta_{Y}+e^{i{\tilde \theta_{k}}}t)}(0)) = 
\zeta_{\Delta_{Y}}(0)+dimker\Delta_{Y}$.

\vskip 0.2 true cm

(2) $\pi_{0}(\Delta_{Y}+\alpha_{k}t)= i\theta_{k}
(\zeta_{\Delta_{Y}}(0)+dimker\Delta_{Y}),$ where $\alpha_{k}=e^{i\theta_{k}}$.

\endproclaim

\vskip 0.2 true cm

Now we are ready to compute $c = \sum_{k}c_{k}$ in Theorem 1.1.
Since $\zeta_{(\Delta_{Y}+\alpha_{k}t)}(0)=\zeta_{(e^{-i\phi_{k}}\Delta_{Y}+e^{i{\tilde \theta_{k}}}t)}(0)$,
from (3.1) and Proposition 3.2 
$$c_{k}=\log2\cdot(\zeta_{\Delta_{Y}}(0)+dimker\Delta_{Y}) + \frac{1}{2} i\theta_{k}(\zeta_{\Delta_{Y}}(0)+dimker\Delta_{Y})$$
and hence
$$\sum_{k}c_{k}=m\log2\cdot(\zeta_{\Delta_{Y}}(0)+dimker\Delta_{Y}).$$
This completes the proof of Theorem 1.2.

\vskip 1 true cm
\S 4 {\bf The adiabatic decomposition of the zeta-determinant of a Laplacian }
\vskip 0.5 true cm
\rm

In this section we are going to prove Theorem 1.4.
Recall that
$$
M_{1,r} =M_{1}\cup_{Y} [-r,0]\times Y, \quad
M_{2,r} =M_{2}\cup_{Y} [0,r]\times Y,
$$
where we identify $Y$ with $\{ -r \}\times Y$ and $Y$ with $\{ r \}\times Y$.
Then
$$M_{r}=M_{1,r}\cup_{\{0\}\times Y} M_{2,r}.  \tag4.1$$
Throughout this section we denote  $\{ r \}\times Y$ by $Y_{r}$,  the Dirichlet
(Neumann) condition on $Y_{r}$ by $B_{r}$ ($C_{r}$) and  the Dirichlet condition on $Z$ by $D$.
We assume that $\Delta_{M,D}$ is invertible. Then under some conditions
$\Delta_{M_{r},D}$ is also invertible for $r$ large enough (Lemma 4.6).

From the decomposition (4.1) and Corollary 1.3, we have
$$ \multline
\log Det\Delta_{M_{r},D}=\log Det\Delta_{M_{1,r},B_{0},D}+\log Det\Delta_{M_{2,r},B_{0},D} \\
-\log2\cdot (\zeta_{\Delta_{Y}}(0)+dimker\Delta_{Y})+\log DetR_{M_{r}}. 
\endmultline \tag4.2
$$
From the decomposition $M_{r}=(M_{1}\cup M_{2})\cup N_{r}$ with $N_{r}=[-r,r]\times Y$, we have
$$\split
\log Det\Delta_{M_{r},D} 
&= \log Det\Delta_{(M_{1}\cup M_{2}),B_{-r},B_{r},D} + \log Det\Delta_{N_{r},B_{-r},B_{r}} \\ 
&\qquad -\log2\cdot (\zeta_{\Delta_{Y \cup Y}}(0)+dimker\Delta_{Y \cup Y})+\log DetR_{-r,r}  \\ 
&=\log Det\Delta_{M_{1},B,D} +\log Det\Delta_{M_{2},B,D} + \log Det\Delta_{N_{r},B_{-r},B_{r}} \\ 
&\qquad -2 \log2\cdot (\zeta_{\Delta_{Y}}(0)+dimker\Delta_{Y})+\log DetR_{-r,r}, 
\endsplit  \tag4.3
$$
where $R_{-r,r}:C^{\infty}(Y_{-r}) \oplus C^{\infty}(Y_{r})
\rightarrow C^{\infty}(Y_{-r}) \oplus C^{\infty}(Y_{r}) $
is the Dirichlet-to-Neumann operator corresponding to the decomposition 
$(M_{1}\cup M_{2})\cup N_{r}$.

Put $N_{-r,0}=[-r,0]\times Y$ and $N_{0,r}=[0,r]\times Y$.
Since 
$M_{1,r}=M_{1}\cup N_{-r,0}$ and $M_{2,r}=M_{2}\cup N_{0,r}$,
we have
$$
\multline
\log Det\Delta_{M_{1,r},B_{0},D}=\log Det\Delta_{M_{1},B,D}+\log Det\Delta_{N_{-r,0},B_{-r},B_{0}} \\
-\log2\cdot (\zeta_{\Delta_{Y}}(0)+dimker\Delta_{Y})+\log DetR_{M_{1,r}}, 
\endmultline \tag4.4
$$
$$
\multline
\log Det\Delta_{M_{2,r},B_{0},D}=\log Det\Delta_{M_{2},B,D}+\log Det\Delta_{N_{0,r},B_{0},B_{r}} \\
-\log2\cdot (\zeta_{\Delta_{Y}}(0)+dimker\Delta_{Y})+\log DetR_{M_{2,r}}.
\endmultline  \tag4.5
$$
Here $\Delta_{N_{-r,0},B_{-r},B_{0}}=-\partial_{u}^{2}+\Delta_{Y}$ with the domain 
$\{\phi\in C^{\infty}(N_{-r,0}) \mid \phi|_{Y_{-r}}=\phi|_{Y_{0}}=0 \}$
and $R_{M_{1,r}}$ is the Dirichlet-to-Neumann operator corresponding to the decomposition
$M_{1,r}=M_{1}\cup ([-r,0]\times Y)$.
$\Delta_{N_{0,r},B_{0},B_{r}}$ and $R_{M_{2,r}}$ are defined similarly.

Then from (4.2) to (4.5), we have
$$
-\log2\cdot (\zeta_{\Delta_{Y}}(0)+dimker\Delta_{Y}) + \log DetR_{M_{r}} = $$
$$ \log Det\Delta_{N_{r},B_{-r},B_{r}} -
\log Det\Delta_{N_{-r,0},B_{-r},B_{0}} 
- \log Det\Delta_{N_{0,r},B_{0},B_{r}} $$
$$ + \log DetR_{-r,r} - \log DetR_{M_{1,r}} - \log DetR_{M_{2,r}}.  \tag4.6
$$
From the decomposition of $N_{r}$ as
$$N_{r}=([-r,0]\times Y)\cup ([0,r]\times Y),$$
we have
$$\multline
\log Det\Delta_{N_{r},B_{-r},B_{r}}-\log Det\Delta_{N_{-r,0},B_{-r},B_{0}}-\log Det\Delta_{N_{0,r},B_{0},B_{r}} \\
=-\log2\cdot (\zeta_{\Delta_{Y}}(0)+dimker\Delta_{Y}) + \log DetR_{N_{r}},
\endmultline \tag4.7
$$
where $R_{N_{r}}:C^{\infty}(Y_{0}) \rightarrow C^{\infty}(Y_{0})$ is defined
as follows.
For $f\in C^{\infty}(Y_{0})$, choose $\phi(u,y)$ so that 
$(-\partial_{u}^{2}+\Delta_{Y})\phi=0$ on $N_{r}-Y_{0}$, $\phi|_{Y_{0}}=f$,
$\phi|_{Y_{-r}}=\phi|_{Y_{r}}=0$.
Then, 
$R_{N_{r}}(f)=\left(\partial_{u}(\phi|_{N_{-r,0}})-\partial_{u}(\phi|_{N_{0,r}})\right)|_{Y_{0}}.$
Hence, we obtain from (4.6) and (4.7)
$$\log DetR_{M_{r}}=\log DetR_{N_{r}}+\log DetR_{-r,r} - \log DetR_{M_{1,r}} - \log DetR_{M_{2,r}}.
\tag4.8 $$

\vskip 0.2 true cm

Now we are going to find the spectrum of $R_{N_{r}}:C^{\infty}(Y_{0}) \rightarrow C^{\infty}(Y_{0})$.
For $f_{k}\in C^{\infty}(Y_{0})$ with $\Delta_{Y}f_{k}=\lambda_{k}f_{k}$, we have
$$
\phi(u,y)=\left\{
\aligned \left(e^{\sqrt{\lambda_{k}}u} + \frac{e^{-\sqrt{\lambda_{k}}r}}
{e^{\sqrt{\lambda_{k}}r}-e^{-\sqrt{\lambda_{k}}r}}(e^{\sqrt{\lambda_{k}}u}-e^{-\sqrt{\lambda_{k}}u})
\right)f_{k}(y) & \text{ for } (u,y)\in N_{-r,0} \\
\left(e^{-\sqrt{\lambda_{k}}u} - \frac{e^{-\sqrt{\lambda_{k}}r}}
{e^{\sqrt{\lambda_{k}}r}-e^{-\sqrt{\lambda_{k}}r}}(e^{\sqrt{\lambda_{k}}u}-e^{-\sqrt{\lambda_{k}}u})
\right)f_{k}(y) & \text{ for } (u,y)\in N_{0,r} .
\endaligned
\right.
$$  
Hence,
$$
R_{N_{r}}(f_{k})=\left( 2\sqrt{\lambda_{k}} + \frac{4\sqrt{\lambda_{k}}e^{-\sqrt{\lambda_{k}}r}}
{e^{\sqrt{\lambda_{k}}r}-e^{-\sqrt{\lambda_{k}}r}} \right)f_{k},
$$   
where we interpret $\frac{4\sqrt{\lambda_{k}}e^{-\sqrt{\lambda_{k}}r}}
{e^{\sqrt{\lambda_{k}}r}-e^{-\sqrt{\lambda_{k}}r}}$ as $\frac{2}{r}$ when $\lambda_{k}=0$.
The spectrum of $R_{N_{r}}$ is
$$\left\{ 2\sqrt{\lambda_{k}} + \frac{4\sqrt{\lambda_{k}}e^{-\sqrt{\lambda_{k}}r}}
{e^{\sqrt{\lambda_{k}}r}-e^{-\sqrt{\lambda_{k}}r}} \mid \lambda_{k}\in Spec(\Delta_{Y}) \right\}.$$  

Let $P_{Ker\Delta_{Y}} : C^{\infty}(Y) \rightarrow C^{\infty}(Y)$ be the orthogonal projection onto $Ker\Delta_{Y}$.
Then
$$\multline
\zeta_{R_{N_{r}}}(s)-\zeta_{(2\sqrt{\Delta_{Y}}+\frac{2}{r}P_{Ker\Delta_{Y}})}(s)= \\
\sum_{\lambda_{k}\neq 0}\left\{ \left( 2\sqrt{\lambda_{k}}+\frac{4\sqrt{\lambda_{k}}e^{-\sqrt{\lambda_{k}}r}}
{e^{\sqrt{\lambda_{k}}r} - e^{-\sqrt{\lambda_{k}}r}} \right)^{-s}- \left( 2\sqrt{\lambda_{k}} \right)^{-s} \right\}.
\endmultline $$
The following lemma can be checked easily.
\proclaim{Lemma 4.1}
Let $A$ be an invertible elliptic operator of order $>0$ and $K_{r}$ be a one-parameter family of 
trace class operators
such that $\lim_{r\to\infty}Tr(K_{r})=0$.
Then
$$\lim_{r\to\infty}\log Det(A+K_{r})=\log DetA.$$
\endproclaim
{\it Proof} 
\hskip 0.3 true cm
Note that
$$\split
\log Det(A+K_{r})-\log DetA
& =\int_{0}^{1}\frac{d}{dt}\log Det(A+tK_{r})dt \\
& =\int_{0}^{1}Tr\left( (A+tK_{r})^{-1} K_{r} \right) dt.   \endsplit $$ 
\noindent
If we denote by $\lambda_{0}$ the smallest eigenvalue of $|A|$, for
$r$ large enough we have
$$| \log Det(A+K_{r})-\log Det A| \leq \frac{1}{2\lambda_{0}} Tr(K_{r})$$
and hence the result follows. 
\qed

\noindent
Applying Lemma 4.1 with $A=2\sqrt{\Delta_{Y}}$ and  
$K_{r}=g_{r}(\Delta_{Y})$ with $g_{r}(x)=
\frac{4\sqrt{x}e^{-\sqrt{x}r}}{e^{\sqrt{x}r} - e^{-\sqrt{x}r}}$
on the orthogonal complement of $Ker\Delta_{Y}$,
we get the following equation.
$$
\lim_{r\to\infty}\left\{ \log DetR_{N_{r}}-\log Det(2\sqrt{\Delta_{Y}}+\frac{2}{r}P_{Ker\Delta_{Y}}) \right\}=0.$$
Since 
$$\multline
\log Det(2\sqrt{\Delta_{Y}}+\frac{2}{r}P_{Ker\Delta_{Y}})=\log2\cdot(\zeta_{\Delta_{Y}}(0)+dimKer\Delta_{Y})
+\frac{1}{2}\log Det\Delta_{Y} \\
-(dimKer\Delta_{Y})\log r, \endmultline \tag4.9 $$
we get the following corollary.

\proclaim {Corollary 4.2}
$$ \multline
\lim_{r\to\infty}\left( \log DetR_{N_{r}}+(dimKer\Delta_{Y})\log r\right) =  \\
\log2\cdot(\zeta_{\Delta_{Y}}(0)+dimKer\Delta_{Y})+\frac{1}{2}\log Det\Delta_{Y}. \endmultline $$
\endproclaim

\vskip 0.2 true cm

Now we discuss the operators $R_{M_{1,r}}$, $R_{M_{2,r}}$ and $R_{-r,r}$.
First, we can describe $R_{M_{1,r}}:C^{\infty}(Y_{-r}) \rightarrow C^{\infty}(Y_{-r})$ as follows. 
For $f_{k}\in C^{\infty}(Y_{-r})$ with $\Delta_{Y}f_{k}=\lambda_{k}f_{k}$, we choose the section
$\phi\in C^{0}(M_{1,r})$ satisfying $\Delta_{M_{1,r}}\phi=0$ on $M_{1,r}-Y_{-r}$,
$\phi|_{Y_{-r}}=f_{k}$ and $\phi|_{Z}=\phi|_{Y_{0}}=0$.
Then one can check that
$$\split
R_{M_{1,r}}(f_{k}) &= Q_{1}(f_{k}) - (\partial_{u}(\phi|_{N_{-r,0}}))|_{Y_{-r}} \\
&= Q_{1}(f_{k}) + \left( \sqrt{\lambda_{k}} + \frac{2\sqrt{\lambda_{k}}e^{-\sqrt{\lambda_{k}}r}}
{e^{\sqrt{\lambda_{k}}r}-e^{-\sqrt{\lambda_{k}}r}} \right)f_{k}.
\endsplit
$$
As the same way
$$\split
R_{M_{2,r}}(f_{k}) &= Q_{2}(f_{k}) + (\partial_{u}(\phi|_{N_{0,r}}))|_{Y_{r}} \\
&= Q_{2}(f_{k}) + \left( \sqrt{\lambda_{k}} + \frac{2\sqrt{\lambda_{k}}e^{-\sqrt{\lambda_{k}}r}}
{e^{\sqrt{\lambda_{k}}r}-e^{-\sqrt{\lambda_{k}}r}} \right)f_{k}.
\endsplit
$$
Similarly, $R_{-r,r}:C^{\infty}(Y_{-r})\oplus C^{\infty}(Y_{r}) 
\rightarrow C^{\infty}(Y_{-r}) \oplus C^{\infty}(Y_{r})$ is described as follows.
$$\multline
R_{-r,r}(f_{k},0) = \\
\left(Q_{1}(f_{k})+ \left( \sqrt{\lambda_{k}} + \frac{2\sqrt{\lambda_{k}}e^{-2\sqrt{\lambda_{k}}r}}
{e^{2\sqrt{\lambda_{k}}r}-e^{-2\sqrt{\lambda_{k}}r}} \right)f_{k}, \hskip 0.2 true cm
-\frac{2\sqrt{\lambda_{k}}}{e^{2\sqrt{\lambda_{k}}r}-e^{-2\sqrt{\lambda_{k}}r}}f_{k}\right),
\endmultline
$$
$$\multline
R_{-r,r}(0,f_{k}) = \\
\left(-\frac{2\sqrt{\lambda_{k}}}{e^{2\sqrt{\lambda_{k}}r}-e^{-2\sqrt{\lambda_{k}}r}}f_{k}, \hskip 0.2 true cm
Q_{2}(f_{k})+ \left( \sqrt{\lambda_{k}} + \frac{2\sqrt{\lambda_{k}}e^{-2\sqrt{\lambda_{k}}r}}
{e^{2\sqrt{\lambda_{k}}r}-e^{-2\sqrt{\lambda_{k}}r}} \right)f_{k}\right).
\endmultline
$$
We, therefore,  have
$$R_{-r,r}=\left( \matrix
Q_{1}+ \sqrt{\Delta_{Y}} & 0 \\ 0 & Q_{2}+ \sqrt{\Delta_{Y}}
\endmatrix \right)
+ h_{r}(\Delta_{Y})
\left( \matrix e^{-2r \sqrt{\Delta_{Y}}} & - 1 \\ - 1 & e^{-2r \sqrt{\Delta_{Y}}}  \endmatrix \right),$$
where $h_{r}(x)= \frac{2 \sqrt{x}}{e^{2r\sqrt{x}}-e^{-2r\sqrt{x}}}$ and 
$h_{r}(\Delta_{Y})$ acts on the $Ker\Delta_{Y}$
as the multiplication of $\frac{1}{2r}$.

We are going to discuss the operators $Q_{i}$ and $Q_{i}+\sqrt{\Delta_{Y}}$. 
The following lemma can be checked by using integration by parts ({\it c.f.} Proposition 4.3 in [2]).

\proclaim{Lemma 4.3}
Suppose that $\nabla$ is a connection which is compatible to the inner product on $M$.
{\it i.e.} for any sections $s_{1}$, $s_{2}\in C^{\infty}(E)$ and a tangent vector $w$,
$w(s_{1},s_{2}) = (\nabla_{w}s_{1},s_{2}) + (s_{1},\nabla_{w}s_{2})$.
If $\Delta_{M}=\nabla^{\ast}\nabla$, then each $Q_{i}$ is a non-negative, self-adjoint operator.
\endproclaim

Next, let us consider a Dirac Laplacian for a Dirac operator $A$ which has the form
$G(\partial_{u}+B)$ near the boundary $Y$, where $G$ is 
a bundle automorphism satisfying the conditions (1.1)
and both $G$ and $B$ do not depend on the normal coordinate $u$.
We refer to [5] for the following lemma ({\it c.f.} Lemma 3.1 in [5]).
\proclaim{Lemma 4.4}
Let $\phi$ and $\psi$ be smooth sections on $M_{j}$ ($j=1$, $2$). Then,
$$\langle A_{M_{j}} \phi, \psi \rangle_{M_{j}} -  \langle  \phi, A_{M_{j}} \psi \rangle_{M_{j}}
= \epsilon_{j} \langle \phi|_{Y}, G(\psi|_{Y}) \rangle_{Y}, $$
where $\epsilon_{j}=1$ for $j=2$ and  $\epsilon_{j}=-1$ for $j=1$.
\endproclaim
\noindent
Suppose that  for $f\in C^{\infty}(Y)$, $\phi_{j}$ is the solution of $A^{2}_{M_{j}}$ with 
$\phi_{j}|_{Y}=f$, $\phi_{j}|_{Z}=0$.
Then by Lemma 4.4 
$$
\langle (Q_{1}+|B|)f, f\rangle_{Y} = \langle A_{M_{1}}\phi_{1}, A_{M_{1}}\phi_{1}\rangle_{M_{1}} + 
\langle (|B|-B)f,f\rangle_{Y},  \tag4.10 $$
$$
\langle (Q_{2}+|B|)f, f\rangle_{Y} = \langle A_{M_{2}}\phi_{2}, A_{M_{2}}\phi_{2}\rangle_{M_{2}} + 
\langle (|B|+B)f,f\rangle_{Y}.   \tag4.11 $$
As a consequence,  $f\in Ker(Q_{1}+|B|)$ if and only if $A_{M_{1}}\phi_{1}=0$ and $f\in Im P_{\geq}$ and hence
on the cylinder part we can express $\phi_{1}$ as
$$\phi_{1}=\sum_{j=1}^{k}a_{j}g_{j}+\sum_{\lambda_{j}>0}b_{j}e^{-\lambda_{j}u}h_{j},$$
where $Bg_{j}=0$, $Bh_{j}=\lambda_{j}h_{j}$.
This implies that $\phi_{1}$ is the restriction of an extended $L^{2}$-solution of $A_{M_{1},\infty}$ on 
$M_{1,\infty}:=M_{1} \cup_{Y} Y \times [0,\infty)$.
We can say similar assertion for $\phi_{2}$ and have the following corollary 
({\it c.f.} Theorem 2.2 in [5], see also [1], [2]).
\proclaim{Corollary 4.5}
The invertibility of $Q_{1}+\sqrt{B^{2}}$ and $Q_{2}+\sqrt{B^{2}}$ is equivalent to the non-existence of the extended $L^{2}$-solutions
of  $A_{M_{1,\infty}}$ and $A_{M_{2,\infty}}$ on $M_{1,\infty}$ and $M_{2,\infty}$. In particular, this condition implies that 
$KerB=0$.
\endproclaim

\vskip 0.3 true cm

\proclaim{Lemma 4.6}
Suppose that $\Delta_{M}$ is either a connection Laplacian or a Dirac Laplacian for a connection
compatible to the inner product as above and $\Delta_{M,D}$ is invertible. 
If both $Q_{1}+\sqrt{\Delta_{Y}}$ and $Q_{2}+\sqrt{\Delta_{Y}}$ are invertible, then 
$R_{-r,r}$ and $\Delta_{M_{r,D}}$ are invertible for $r$ large enough.
\endproclaim
{\it Proof} \hskip 0.3 true cm
We are going to show first that $R_{-r,r}$ is injective. Then this implies that $\Delta_{M_{r,D}}$ is injective. 
Since $\Delta_{M_{r,D}}$ is self-adjoint, $\Delta_{M_{r,D}}$ is invertible and this implies again that
$R_{-r,r}$ is also invertible ([3], [8]).  

Putting $A_{r}=h_{r}(\Delta_{Y})$ with
$h_{r}(x)=\frac{2 \sqrt{x}}{e^{2r\sqrt{x}}-e^{-2r\sqrt{x}}}$,
$$
\left< R_{-r,r}\binom{f}{g}, \binom{f}{g}\right>_{L^{2}(Y)} $$
$$ \multline 
=\langle (Q_{1}+\sqrt{\Delta_{Y}})f,f\rangle + \langle (Q_{2}+\sqrt{\Delta_{Y}})g,g\rangle  \\
+ \langle A_{r}e^{-2r\sqrt{\Delta_{Y}}}f,f\rangle + \langle A_{r}e^{-2r\sqrt{\Delta_{Y}}}g,g\rangle
- \langle A_{r}g,f\rangle- \langle A_{r}f,g\rangle. \endmultline$$
Note that each $Q_{i}+\sqrt{\Delta_{Y}}$ is a non-negative operator by Lemma 4.3 and (4.10), (4.11).
Let $\lambda_{0}$ be the minimum of the eigenvalues of $Q_{1}+\sqrt{\Delta_{Y}}$ and $Q_{2}+\sqrt{\Delta_{Y}}$.
Since $\lim_{r\to\infty}||A_{r}||_{L^{2}}=0$, one can choose $r_{0}$ so that for $r\geq r_{0}$,  $||A_{r}||_{L^{2}} < \lambda_{0}$. 
Then for $r\geq r_{0}$, $R_{-r,r}$ is injective and this completes the proof.
\qed

\vskip 0.2 true cm
In case that both $Q_{1}+\sqrt{\Delta_{Y}}$ and $Q_{2}+\sqrt{\Delta_{Y}}$ are invertible, we can apply Lemma 4.1 directly.
\proclaim{Corollary 4.7} 
\hskip 0.2 true cm
Assume that  both $Q_{1}+\sqrt{\Delta_{Y}}$ and $Q_{2}+\sqrt{\Delta_{Y}}$ are invertible. Then :

\vskip 0.2 true cm
$ (1) \lim_{r\to\infty}\log DetR_{M_{1,r}} = \log Det(Q_{1}+\sqrt{\Delta_{Y}}).$
\vskip 0.2 true cm
$ (2) \lim_{r\to\infty}\log DetR_{M_{2,r}}= \log Det(Q_{2}+\sqrt{\Delta_{Y}}).$
\vskip 0.2 true cm
$(3) \lim_{r\to\infty}\log DetR_{-r,r}=\log Det(Q_{1}+\sqrt{\Delta_{Y}})+\log Det(Q_{2}+\sqrt{\Delta_{Y}}).$
\endproclaim

\vskip 0.2 true cm

Combining Corollary 4.2 with Corollary 4.7 and (4.2), (4.8), we complete the proof of Theorem 1.4.

\vskip 1 true cm

\S 5 {\bf The adiabatic decomposition of the analytic torsion}

\vskip 0.5 true cm

In this section, we are going to prove Theorem 1.6 and Theorem 1.7.
Recall that $M$ is a closed manifold of dimension $m$ with the product structure near a hypersurface $Y$.
We define $M_{r}$, $M_{1,r}$ and $M_{2,r}$ as in Section 4 and
suppose that $\rho_{M_{r}}$ ($\rho_{M_{1,r}}$, $\rho_{M_{2,r}}$, $\rho_{Y}$) is an orthogonal representation
of $\pi_{1}(M_{r})$ ($\pi_{1}(M_{1,r})$, $\pi_{1}(M_{2,r})$, $\pi_{1}(Y)$) to $SO(n)$, respectively.
Then we can construct a flat bundle $E_{\rho_{M_{r}}}={\tilde M_{r}}\times_{\rho_{M_{r}}} {\Bbb R}^{n}$,
where ${\tilde M_{r}}$ is the universal cover of $M_{r}$. The flat bundles $E_{\rho_{M_{1,r}}}$, 
$E_{\rho_{M_{2,r}}}$ and $E_{\rho_{Y}}$ are defined in the same way.

For each $q$, denote by $\Delta_{M_{r}}^{q} :=(d_{q}+d_{q}^{\ast})^{2}$ 
the Hodge Laplacian acting on $q$-forms valued in $E_{\rho_{M_{r}}}$.
Then the analytic torsion $\tau(M_{r},\rho_{M_{r}})$ is defined by
$$\log\tau(M_{r},\rho_{M_{r}}) =\frac{1}{2}\sum_{q=0}^{m}(-1)^{q} \cdot q\cdot \log Det\Delta_{M_{r}}^{q}.$$
\noindent
To define the analytic torsion on $M_{i,r}$, we choose
the absolute or the relative boundary condition on $Y_{0}$. 
Near $Y_{0}$, a differential $q$-form $\omega$ can be expressed by
$$\omega=\omega_{1}+du\wedge\omega_{2}, \tag5.1$$
where $\omega_{1}$ and $\omega_{2}$ do not contain $du$.

\proclaim{Definition 5.1} Suppose that a $q$-form $\omega$ in $M_{1,r}$ is expressed as in (5.1).\newline
(1) $\omega$ satisfies the absolute boundary condition 
if $(\partial_{u}\omega_{1})|_{Y_{0}}=0$ and $\omega_{2}|_{Y_{0}}=0$. \newline
(2) $\omega$ satisfies the relative boundary condition if  $\omega_{1}|_{Y_{0}}=0$ and 
$(\partial_{u}\omega_{2})|_{Y_{0}}=0$. 
\endproclaim

We denote by $\Omega_{abs}^{q}(M_{i,r})$, $\Omega_{rel}^{q}(M_{i,r})$
the sets of all $q$-forms valued in $E_{M_{i,r}}$ satisfying the absolute and the relative boundary
conditions, respectively.
We also denote by $\Delta_{M_{i,r},abs}^{q}$, 
$\Delta_{M_{i,r},rel}^{q}$ the Laplacian acting on $q$-forms 
valued in $E_{M_{i,r}}$ with
$$Dom(\Delta_{M_{i,r},abs}^{q})=\Omega_{abs}^{q}(M_{i,r}), \qquad\qquad 
Dom(\Delta_{M_{i,r},rel}^{q})=\Omega_{rel}^{q}(M_{i,r}).$$
Then the analytic torsion $\tau_{abs}(M_{i,r},\rho_{M_{i,r}})$ and 
$\tau_{rel}(M_{i,r},\rho_{M_{i,r}})$ are defined by
$$\log\tau_{abs}(M_{i,r},\rho_{M_{i,r}})= \frac{1}{2}\sum_{q=0}^{m}(-1)^{q}\cdot q\cdot 
\log Det\Delta_{M_{i,r},abs}^{q},$$
$$\log\tau_{rel}(M_{i,r},\rho_{M_{i,r}})= \frac{1}{2}\sum_{q=0}^{m}(-1)^{q}\cdot q\cdot 
\log Det\Delta_{M_{i,r},rel}^{q}.$$
It is a well-known fact ({\it c.f.} [11]) that 
$$Ker\Delta^{q}_{M_{i,r},abs} \cong H^{q}(M_{i,r};\rho_{M_{i,r}}) ,\qquad 
Ker\Delta^{q}_{M_{i,r},rel} \cong H^{q}(M_{i,r},Y;\rho_{M_{i,r}}).$$

We consider $M_{1,r}$ (a manifold with boundary $Y$), first. 
Recall that $M_{1,r}=M_{1}\cup_{Y_{-r}} N_{-r,0}$ with $N_{-r,0}=[-r,0]\times Y$, and $Y_{-r}= \{ -r \} \times Y$,
$Y_{0}= \{ 0 \} \times Y$.
We denote by $B$, $D$, the Dirichlet boundary conditions on $Y_{-r}$, $Y_{0}$, respectively.

For a given representation $\rho_{M_{1,r}} : \pi_{1}(M_{1,r})\rightarrow SO(n)$, define
$\rho_{Y}:\pi_{1}(Y)\rightarrow SO(n)$ by $\rho_{Y}=\rho_{M_{1,r}}\circ \iota_{Y}$, where 
$\iota_{Y}:\pi_{1}(Y)\rightarrow \pi_{1}(M_{1,r})$ is the natural homomorphism. 
Then the restriction of the bundle $E_{\rho_{M_{1,r}}}$ to $Y$ is isomorphic to $E_{\rho_{Y}}$, 
({\it c.f.} [11]).

The set $\Omega^{q}(N_{-r,0},E_{\rho_{M_{1,r}}}|_{N_{-r,0}})$
of $q$-forms valued in $E_{\rho_{M_{1,r}}}|_{N_{-r,0}}$ can be decomposed as follows.
$$\multline
\Omega^{q}(N_{-r,0},E_{\rho_{M_{1,r}}}|_{N_{-r,0}})=
C^{\infty}([-r,0],E_{\rho_{M_{1,r}}}|_{N_{-r,0}})\otimes \Omega^{q}(Y,E_{\rho_{Y}}) \oplus \\
du\wedge C^{\infty}([-r,0],E_{\rho_{M_{1,r}}}|_{N_{-r,0}})\otimes \Omega^{q-1}(Y,E_{\rho_{Y}}).
\endmultline \tag5.2 $$
From this decomposition, the Laplacian $\Delta^{q}_{M_{1,r}}$, when restricted to $N_{-r,0}$, can be
expressed by
$$\Delta^{q}_{M_{1,r}} = -\partial_{u}^{2} + \left(\matrix \Delta_{Y}^{q} & 0 \\
0 & \Delta_{Y}^{q-1} \endmatrix \right),  \tag5.3$$
where $\Delta_{Y}^{q}$ is the Laplacian acting on $q$-forms on $Y$, valued in $E_{\rho_{Y}}$. Here
and throughout this section we use the convention that $\Delta_{Y}^{q}=0$ for $q<0$ or $q\geq m$. 

To describe the gluing formula of 
the type of Theorem 1.1 (or Corollary 1.3) in this context, we need to define 
modified Dirichlet-to-Neumann operators $Q_{1}^{q}$, $Q_{N_{-r,0},abs}^{q}$ and $Q_{N_{-r,0},rel}^{q}$ as follows. For  simplicity denote 
$E=(\wedge^{q}T^{\ast}M_{1,r}) \otimes E_{\rho_{M_{1,r}}}$.
For a given $f\in C^{\infty}(E|_{Y_{-r}})$, choose smooth sections $\phi\in C^{\infty}(E|_{M_{1}})$,
$\psi_{abs}\in C^{\infty}(E|_{N_{-r,0}})$ and $\psi_{rel}\in C^{\infty}(E|_{N_{-r,0}})$ such that
$$\Delta^{q}_{M_{1}}\phi=0, \quad \Delta^{q}_{N_{-r,0}}\psi_{abs}=\Delta^{q}_{N_{-r,0}}\psi_{rel}=0, \quad
\phi|_{Y_{-r}}=\psi_{abs}|_{Y_{-r}}=\psi_{rel}|_{Y_{-r}}=f,$$                             
and
$\psi_{abs}$ ($\psi_{rel}$) satisfies the absolute (relative) boundary condition on $Y_{0}$, respectively. Then we define
$$Q_{1}^{q}(f)=(\partial_{u}\phi)|_{Y_{-r}},$$ 
$$Q^{q}_{N_{-r,0},abs}(f)=(-\partial_{u}\psi_{abs})|_{Y_{-r}}, \quad
Q^{q}_{N_{-r,0},rel}(f)=(-\partial_{u}\psi_{rel})|_{Y_{-r}},$$
and
$$R^{q}_{B,abs}=Q_{1}^{q}+Q^{q}_{N_{-r,0},abs}, \quad R^{q}_{B,rel}=Q_{1}^{q}+Q^{q}_{N_{-r,0},rel}.$$
Then the following theorem can be proved in the same way as Theorem 1.1 ({\it c.f.} see the {\it Remark}
after Corollary 1.3).

\proclaim{Theorem 5.2}  We denote $k_{q}=dimKer\Delta_{Y}^{q}$. Then :\newline
$$\multline
(1)\quad \log Det\Delta_{M_{1,r},D}^{q}- \log Det\Delta_{M_{1},B}^{q}-\log Det\Delta_{N_{-r,0},B,D}^{q} \\
=-\log2(\zeta_{\Delta_{Y}^{q-1}}(0) +\zeta_{\Delta_{Y}^{q}}(0)+k_{q-1}+k_{q}) + \log DetR^{q}_{B,D}.
\endmultline $$
$$\multline
(2)\quad \log Det\Delta_{M_{1,r},abs}^{q}- \log Det\Delta_{M_{1},B}^{q}-\log Det\Delta_{N_{-r,0},B,abs}^{q} \\
=-\log2(\zeta_{\Delta_{Y}^{q-1}}(0)+\zeta_{\Delta_{Y}^{q}}(0)+k_{q-1}+k_{q}) + \log DetR^{q}_{B,abs}.
\endmultline $$
$$\multline
(3)\quad \log Det\Delta_{M_{1,r},rel}^{q}- \log Det\Delta_{M_{1},B}^{q}-\log Det\Delta_{N_{-r,0},B,rel}^{q} \\
=-\log2(\zeta_{\Delta_{Y}^{q-1}}(0)+\zeta_{\Delta_{Y}^{q}}(0)+k_{q-1}+k_{q}) + \log DetR^{q}_{B,rel}.
\endmultline $$
\endproclaim

We next describe the operators $\Delta_{N_{-r,0},B,abs}^{q}$ and $\Delta_{N_{-r,0},B,rel}^{q}$.
From the decomposition (5.2), we have :
$$
\Delta^{q}_{N_{-r,0},B,abs}= \left( \matrix (-\partial_{u}^{2}+\Delta_{Y}^{q})_{N_{-r,0},B,C} & 0 \\
0 & (-\partial_{u}^{2}+\Delta_{Y}^{q-1})_{N_{-r,0},B,D} \endmatrix \right), $$
$$
\Delta^{q}_{N_{-r,0},B,rel}= \left( \matrix (-\partial_{u}^{2}+\Delta_{Y}^{q})_{N_{-r,0},B,D} & 0 \\
0 & (-\partial_{u}^{2}+\Delta_{Y}^{q-1})_{N_{-r,0},B,C} \endmatrix \right), $$
where $C$ means the Neumann boundary condition on $Y_{0}$ and $B$ ($D$) means the Dirichlet boundary
condition on $Y_{-r}$ ($Y_{0}$).
Hence, we have :

$$ \multline
\log Det\Delta^{q}_{N_{-r,0},B,abs}-\log Det\Delta^{q}_{N_{-r,0},B,D}= \\
\log Det(-\partial_{u}^{2}+\Delta_{Y}^{q})_{N_{-r,0},B,C}-
\log Det(-\partial_{u}^{2}+\Delta_{Y}^{q})_{N_{-r,0},B,D}, \endmultline  \tag5.4 $$
$$ \multline
\log Det\Delta^{q}_{N_{-r,0},B,rel}-\log Det\Delta^{q}_{N_{-r,0},B,D}= \\
\log Det(-\partial_{u}^{2}+\Delta_{Y}^{q-1})_{N_{-r,0},B,C}-
\log Det(-\partial_{u}^{2}+\Delta_{Y}^{q-1})_{N_{-r,0},B,D}. \endmultline \tag5.5
$$
Now we assume that 
$Q_{1}^{q} +\left( \smallmatrix \sqrt{\Delta_{Y}^{q}} & 0 \\  0 & \sqrt{\Delta_{Y}^{q-1}} \endsmallmatrix \right)$ is invertible. 
Since the Hodge Laplacian $\Delta_{M_{r}}^{q}$ is a Dirac Laplacian satisfying (1.1), 
by Corollary 4.5 we have $k_{q-1}=k_{q}=0$ ({\it c.f.} see the {\it Remark} below Theorem 1.6). 
By Corollary 1.5 and (5.4), (5.5) we have :
$$
\lim_{r\to\infty}\left\{ \log Det\Delta^{q}_{N_{-r,0},B,abs}-\log Det\Delta^{q}_{N_{-r,0},B,D} 
\right\} = \frac{1}{2}\log Det\Delta_{Y}^{q}, \tag5.6$$
$$
\lim_{r\to\infty}\left\{ \log Det\Delta^{q}_{N_{-r,0},B,rel}-\log Det\Delta^{q}_{N_{-r,0},B,D} 
\right\} = \frac{1}{2}\log Det\Delta_{Y}^{q-1}. \tag5.7$$
From Theorem 5.2 , we have :
$$ \log Det\Delta^{q}_{M_{1,r},abs} - \log Det\Delta^{q}_{M_{1,r},D} $$
$$\multline
= \left( \log Det\Delta^{q}_{N_{-r,0},B,abs} - \log Det\Delta^{q}_{N_{-r,0},B,D} \right) \\
+ \log DetR^{q}_{B,abs} - \log DetR^{q}_{B,D}, \endmultline \tag5.8$$
and
$$ \log Det\Delta^{q}_{M_{1,r},rel} - \log Det\Delta^{q}_{M_{1,r},D}  $$
$$ \multline
= \left( \log Det\Delta^{q}_{N_{-r,0},B,rel} - \log Det\Delta^{q}_{N_{-r,0},B,D}  \right) \\
+ \log DetR^{q}_{B,rel} - \log DetR^{q}_{B,D}. \endmultline \tag5.9$$

\proclaim{Lemma 5.3}
Suppose that $Q_{1}^{q} +\left( \smallmatrix \sqrt{\Delta_{Y}^{q}} & 0 \\  0 & \sqrt{\Delta_{Y}^{q-1}} \endsmallmatrix \right)$ 
is invertible. Then,
$$ 
\lim_{r\to\infty}\log Det R^{q}_{B,abs}= \lim_{r\to\infty}\log Det R^{q}_{B,rel}=
\lim_{r\to\infty}\log Det R^{q}_{B,D} $$
$$= \log Det(Q_{1}^{q} +\left( \smallmatrix \sqrt{\Delta_{Y}^{q}} & 0 \\  0 & \sqrt{\Delta_{Y}^{q-1}} \endsmallmatrix \right)).$$
\endproclaim
{\it Proof} \hskip 0.3 true cm
The last equality is exactly the assertion (1) in Corollary 4.7. We are going to show that
$\lim_{r\to\infty}\log Det R^{q}_{B,abs}= 
\log Det\left(Q_{1}^{q} +\left( \smallmatrix \sqrt{\Delta_{Y}^{q}} & 0 \\  0 & \sqrt{\Delta_{Y}^{q-1}} \endsmallmatrix \right)\right)$. 
The case of $\log Det R^{q}_{B,rel}$ can be proved in the same way.

By a direct computation one can check the followings.
For $f\in \Omega^{q}(Y,E_{\rho_{Y}})$ with $\Delta_{Y}^{q}f=\lambda f$,
$$
 R^{q}_{B,abs}(f)=Q_{1}^{q}(f)+\left(\sqrt{\lambda}-\frac{2\sqrt{\lambda}e^{-\sqrt{\lambda}r}}
{e^{\sqrt{\lambda}r}+e^{-\sqrt{\lambda}r}}\right)f.$$
For $g\in \Omega^{q-1}(Y,E_{\rho_{Y}})$ with $\Delta_{Y}^{q-1}g=\mu g$,
$$
 R^{q}_{B,abs}(du\wedge g)=Q_{1}^{q}(du\wedge g)+\left(\sqrt{\mu}+\frac{2\sqrt{\mu}e^{-\sqrt{\mu}r}}
{e^{\sqrt{\mu}r}-e^{-\sqrt{\mu}r}}\right)du\wedge g.$$
Then the result follows from Lemma 4.1.
\qed  \newline
\noindent
From (5.6) to (5.9) and Lemma 5.3, we have the following corollary.
\proclaim{Corollary 5.4}
Suppose that $Q_{1}^{q} +\left( \smallmatrix \sqrt{\Delta_{Y}^{q}} & 0 \\  0 & \sqrt{\Delta_{Y}^{q-1}} \endsmallmatrix \right)$ 
is invertible for each $q$.  Then the following equalities hold.
$$
(1) \lim_{r\to\infty}\left\{ \log Det\Delta^{q}_{M_{1,r},abs}-\log Det\Delta^{q}_{M_{1,r},D} \right\}= 
\left\{ \aligned
\frac{1}{2}\log Det\Delta^{q}_{Y} &  \quad (0\leq q \leq m-1) \\
0  &  \quad (q=m) .\endaligned \right.  $$
$$
(2) \lim_{r\to\infty}\left\{ \log Det\Delta^{q}_{M_{1,r},rel}-\log Det\Delta^{q}_{M_{1,r},D} \right\}=
\left\{
\aligned \frac{1}{2}\log Det\Delta^{q-1}_{Y}  &  \quad  (1 \leq q \leq m) \\ 
0  &  \quad (q=0) . \endaligned \right.    $$
\endproclaim

Now we are ready to prove Theorem 1.6.
$$
\lim_{r\to\infty}\left\{ \log\tau_{abs}(M_{1,r},\rho_{M_{1,r}}) - \log\tau_{rel}(M_{1,r},\rho_{M_{1,r}}) \right\}
$$
$$=\lim_{r\to\infty}
\frac{1}{2}\sum_{q=0}^{m}(-1)^{q}\cdot q \cdot \left(\log Det\Delta_{M_{1,r},abs}-\log Det\Delta_{M_{1,r},D}
\right) \qquad \qquad \qquad
$$
$$\qquad  -\lim_{r\to\infty}
\frac{1}{2}\sum_{q=0}^{m}(-1)^{q}\cdot q \cdot \left(\log Det\Delta_{M_{1,r},rel}-\log Det\Delta_{M_{1,r},D}
 \right)
$$
$$
=\frac{1}{4}\sum_{q=0}^{m-1}(-1)^{q}\cdot q \cdot \log Det(\Delta^{q}_{Y})
- \frac{1}{4}\sum_{q=1}^{m}(-1)^{q}\cdot q \cdot \log Det(\Delta^{q-1}_{Y})  \qquad \qquad \qquad \qquad
$$
$$=\frac{1}{2}\sum_{q=0}^{m-1}(-1)^{q}\cdot q \cdot \log Det(\Delta^{q}_{Y})
+\frac{1}{4}\sum_{q=0}^{m-1}(-1)^{q}\cdot \log Det(\Delta_{Y}^{q}) \qquad \qquad \qquad  \qquad  \qquad $$
$$=\tau(Y,\rho_{Y}). \qquad \qquad \qquad \qquad \qquad \qquad \qquad \qquad \qquad \qquad \qquad \qquad \qquad 
\qquad \qquad $$
This completes the proof of Theorem 1.6.
\qed

\vskip 0.3 true cm

Next, we take care of the closed manifold $M_{r}=M_{1,r}\cup_{Y_{0}} M_{2,r}$.
From Theorem 1.4, we have :
$$\lim_{r\to\infty}\left\{ \log Det\Delta_{M_{r}}^{q}-\log Det\Delta_{M_{1,r},D}^{q}
-\log Det\Delta_{M_{2,r},D}^{q}  \right\} $$
$$=\frac{1}{2}\left(\log Det\Delta^{q}_{Y}+\log Det\Delta^{q-1}_{Y}\right).$$
On the other hand,
$$\lim_{r\to\infty}\left\{ \log Det\Delta_{M_{r}}^{q}-\log Det\Delta_{M_{1,r},D}^{q}
-\log Det\Delta_{M_{2,r},D}^{q}  \right\} $$
$$
= \lim_{r\to\infty} \left\{ \left( \log Det\Delta_{M_{r}}^{q}-\log Det\Delta_{M_{1,r},abs}^{q}-
\log Det\Delta_{M_{2,r},rel}^{q} \right) \right.
$$
$$
+\left(\log Det\Delta_{M_{1,r},abs}^{q}-\log Det\Delta_{M_{1,r},D}^{q}  \right) 
$$
$$
+\left.\left( \log Det\Delta_{M_{2,r},rel}^{q}-\log Det\Delta_{M_{2,r},D}^{q}  \right)\right\}.
$$
From Corollary 5.4, we have :
$$\lim_{r\to\infty}\left\{ \log Det\Delta_{M_{r}}^{q}-\log Det\Delta_{M_{1,r},abs}^{q}-
\log Det\Delta_{M_{2,r},rel}^{q} \right\}=0
$$
and therefore, we obtain
$$\lim_{r\to\infty}\left\{ \tau(M_{r},\rho_{M_{r}}) - \tau_{abs}(M_{1,r},\rho_{M_{1,r}})
-\tau_{rel}(M_{2,r},\rho_{M_{2,r}}) \right\} = 0,$$
which completes the proof of Theorem 1.7.

\vskip 1 true cm

\vskip 0.3 true cm
E-mail address : ywlee\@math.inha.ac.kr

\enddocument